\documentclass[11pt]{article}
\usepackage{graphicx}

\def\la{\big\langle}
\def\ra{\big\rangle}

\def\ds{\displaystyle}

\def\forall{\hbox{for all}~}
\def\I{{\mathcal I}}
\def\H{{\bf H}}

\def\bfv{{\bf v}}
\def\bfw{{\bf w}}
\def\bfz{{\bf z}}
\def\bfu{{\bf u}}
\def\bfq{{\bf q}}

\def\bfh{{\bf h}}
\def\bfk{{\bf k}}
\def\bfg{{\bf g}}
\def\ve{\varepsilon}

\def\Q{{\cal Q}}

\def\R{I\!\!R}
\def\bega{\begin{array}}
\def\enda{\end{array}}

\def\bfe{{\bf e}}

\def\P{{\cal P}}
\def\O{{\cal O}}
\def\C{{\cal C}}
\def\V{{\bf V}}

\def\Tilde{\widetilde}

\def\bfp{{\bf p}}

\def\bel{\begin{equation}\label}
\def\eeq{\end{equation}}

\def\ve{\varepsilon}

\def\T{{\mathcal T}}
\def\I{{\mathcal I}}

\def\Q{{\mathcal Q}}
\def\bfR{{\bf R}}
\def\C{{\mathcal C}}

\def\U{{\mathcal U}}
\def\bfx{{\bf x}}
\def\g{{\bf g}}

\def\p{{\bf p}}
\def\q{{\bf q}}
\def\w{{\bf w}}

\def\H{{\mathcal H}}

\def\sqr#1#2{\vbox{\hrule height .#2pt
\hbox{\vrule width .#2pt height #1pt \kern #1pt
\vrule width .#2pt}\hrule height .#2pt }}
\def\square{\sqr74}
\def\endproof{\hphantom{MM}\hfill\llap{$\square$}\goodbreak}

\def\U{{\cal U}}
\def\T{{\cal T}}
\def\O{{\cal O}}
\def\la{\Big\langle}
\def\ra{\Big\rangle}

\def\C{{\mathcal C}}
\def\R{I\!\!R}

\def\rr{I\!\!R}

\def\bv{{\bf v}}

\def\P{{\mathcal P}}
\def\vs{\vskip 2em}

\def\v{\vskip 1em}

\def\begi{\begin{itemize}}
\def\endi{\end{itemize}}
\def\Tilde{\widetilde}
\def\forall{\hbox{for all }~}
\def\be{\begin{equation}}
\def\beq{\begin{equation}}
\def\bel{\begin{equation}\label}
\def\eeq{\end{equation}}
\def\om{\omega}
\def\Om{{\bf\Omega}}

\newtheorem{thm}{Theorem}[section]

\newtheorem{lma}{Lemma}[section]
\newtheorem{prop}{Proposition}[section]
\newtheorem{remark}{Remark}[section]
\newtheorem{definition}{Definition}[section]

\title{\bf On the Control of Non Holonomic Systems\\
by Active Constraints}
\begin{document}

\vs
\author{Alberto Bressan$^{(*)}$, Ke Han$^{(*)}$,
and Franco Rampazzo$^{(**)}$\\
\\
(*)~Department of Mathematics, Penn State University,
\\University Park, Pa.~16802, USA.\\
(**)  Dipartimento di Matematica Pura ed Applicata,
Universit\`a di Padova,\\ Padova  35141, Italy.
\\E-mails: bressan@math.psu.edu~,~~kxh323@psu.edu~,~ ~rampazzo@math.unipd.it\\}
\maketitle
\vs

\begin{abstract}
The paper is concerned with mechanical systems
which are controlled by implementing a number of time-dependent,
frictionless holonomic constraints.
The main novelty is due to the presence
of additional non-holonomic constraints.
We develop a general framework to analyze these problems,
deriving the equations of motion and studying the
continuity properties of the ``control-to-trajectory" maps.
Various geometric characterizations are provided, in order
that the equations be affine w.r.t.~the time derivative of the control.
In this case the system is {\it fit for jumps}, and the evolution is well defined also in connection with discontinuous control functions.
The classical Roller Racer provides
an example where the non-affine dependence of the equations
on the derivative of the control is due only to the non-holonomic constraint.
This is a case where the presence of quadratic terms in the equations can be used
for controllability purposes.
\end{abstract}

\section{Introduction}
\setcounter{equation}{0}

The control of  mechanical systems provides a rich area
for mathematical investigation.
In a commonly adopted model
\cite{Bloch, BMZ, BL, NS},
the controller can modify the time evolution of
the system by applying
additional external  forces.
This  leads to a control problem in standard form,
where vector fields governing the state variables depend continuously
on the control function.

In an alternative model, also physically meaningful,
the controller acts on the system by directly assigning the
values of some of the coordinates, as functions of time.
Stated in a more intrinsic fashion, this means that the controller
assigns the leaves of a foliation of the configuration manifold as functions of time.
Here the basic framework consists of a manifold with
coordinates $q= (q^1,\ldots, q^N, q^{N+1}, \ldots, q^{N+M})$.
We assume that the values of the last $M$ coordinates can be prescribed:
$q^{N+1} = u^1(t)$, $\ldots\,$, $q^{N+M} = u^M(t)$.
Instead of forces, these controls thus take the form of
time-dependent holonomic constraints.
The evolution of the remaining free coordinates
$q^1,\ldots, q^N$ is then

determined by a control system where
the right hand side depends not only on the control itself, but also
(linearly or quadratically)
on the time derivative of the control
function.

This alternative point of view
was introduced, independently,  in \cite{AB1} and in \cite{Marle}.
A considerable amount of literature is now available on
mechanical systems controlled by active constraints.
The form of the basic equations of motion, in relation with geometrical
properties of the system, was studied in \cite{Baillieul,
AB1, CF, LR, Rampazzo1, Rampazzo2}.
Differential equations of impulsive nature, where the
right hand side contains
a measure, such as the distributional derivative of a
discontinuous control,
were considered in \cite{B-R1, B-R2, B-R3, LS, miller, Rampazzo3, Su1}.
These works focus, in particular, on
the continuity of the ``control-to-trajectory" map,
in various topologies.
In addition, problems of stabilization and of optimal control for
this kind of ``hyper-impulsive" mechanical systems
were studied in \cite{Baillieul, B-R4, L2} and in
\cite{B-R2, B-R3, AB2, motta},
respectively.
  See also \cite{B1} for a survey.

The goal of the present paper is to extend this theory
to the case where
the system is subject to some additional non-holonomic constraints.
In this case, the
velocity vector $\dot q= (\dot q^1,\ldots, \dot q^{N+M})$
satisfies
an additional set of $\nu$ linear relations
$$\sum_{i=1}^{N+M}~ \omega^k_i\, \dot q^i ~=~0
\qquad\qquad k= 1,\ldots,\nu.$$



A major focus of our analysis is on the form of the
resulting equations,
In general, the right hand side of the
evolution equations turns out to be
a {\em quadratic} polynomial of the time derivative
$\dot u(\cdot)$.  Therefore, the
evolution problem is well posed as soon as the
control function satisfies $u(\cdot)\in W^{1,2}$.
There are, however, important cases where the right hand side is
an {\em affine} function of $\dot u(\cdot)$.
Our main result in this direction,
Theorem~\ref{udotlinear}, yields a number of
equivalent analytic and geometric conditions for this to happen.
In the positive case, the ``control-to-trajectory" map can be
extended by continuity to a larger family of
(possibly discontinuous) control functions.
Following \cite{AB1}, we then say that the system is ``fit for jumps".
Section~7 provides an additional geometric characterization
of this property.  Let
\bel{leaf}\Lambda_u~\doteq~\{(q^1,\ldots, q^{N+M})\,;~~q^{N+\alpha} = u^\alpha\,,
~1\leq\alpha\leq M\}\eeq
be the leaf of the foliation corresponding to the control value $u$.
Roughly speaking, we show that the system is ``fit for jumps" if and only if
the ``infinitesimal non-holonomic distance" between leaves of the foliation
remains constant, along directions compatible with the non-holonomic constraint.

For sake of illustration, in the last section we show how the
present  framework applies to the control of the Roller Racer.
An  interesting feature of this example is that,
without the non-holonomic constraint, the equations of motion would be
affine w.r.t.~the time derivative $\dot u$ of the control function.
However, the additional non-holonomic constraint
renders the system {\it not} fit for jumps. It is indeed the presence of a
quadratic term in the derivative of the control that makes forward motion possible.

\section{Non holonomic systems with active constraints as controls}
\setcounter{equation}{0}

Let $N, M, \nu$ be positive integers such that
$\nu\leq N+M$.
Let $\Q$ be an $(N+M)$-dimensional differential manifold,
which will be  regarded as the space of  configurations of
a mechanical system.

Let $\Gamma$ be a distribution on $\Q$, i.e.~a vector sub-bundle
of the tangent space $T\Q$. Throughout the following, we
consider trajectories
$t\mapsto \bfq(t)\in \Q$ of the mechanical system
which are continuously differentiable  and
satisfy the geometric constraint
 \bel{nhol}
  \dot\bfq(t) ~\in ~\Gamma_{\bfq(t)}\,.
  \eeq
  We do not assume  $\Gamma$ to be  integrable,
so that in general (\ref{nhol}) is a {\it non holonomic constraint}.

In our model, the system will be controlled by means of an active
(holonomic, time-dependent) constraint. To describe
this constraint, let an $M$-dimensional differential manifold
  $\U$ be given, together with  a submersion
  \bel{submersion}
  \pi : \Q\mapsto \U.
  \eeq
  The fibers $\pi^{-1}(\bfu)\subset\Q$ will be regarded
  as the states of the active constraint. The set
  of all
  these fibers can be identified with the {\it control manifold}
  $\U$.

  Let $\I$ be a time interval and let $\bfu: \I \mapsto
  \U$ be a continuously differentiable map.
  We say that {\it a trajectory $\bfq:\I\to\Q$ agrees with
  the control $\bfu(\cdot)$ if
  \bel{hol}
  \pi\circ\bfq(t)~ = ~\bfu(t)\qquad\qquad \forall  t\in \I\,.
  \eeq
 }
For each $\bfq\in\Q$, consider the subspace
of the tangent space at $\bfq$ given by
$$\Delta_{\bfq}~\doteq~ ker \,\,T_{\bfq}\pi\,.
$$
Here $T_\bfq\pi$ denotes the linear tangent map between the tangent
spaces $T_\bfq\Q$ and $T_{\pi(\bfq)}\U$.
Clearly, $\Delta$ is the (holonomic) distribution whose integral
manifolds are precisely the fibers $\pi^{-1}(\bfu)$.

\v

 \subsection{General setting}
  \begin{itemize}

   \item[{\bf 1)}] The manifold
  $\Q$ is endowed with a Riemannian
  metric $\g=\g_{\bfq}[\cdot,\cdot]$, the so-called
  {\it kinetic metric},
which defines the kinetic energy $\T$. More precisely,
for each $\bfq\in \Q$ and $\bfv\in T_{\bfq}\Q$
one has
\bel{kinetic}
\T(\bfq,\bfv)~\doteq~ \frac{1}{2} \g_{\bfq}[\bfv,\bfv]\,.
\eeq
We shall use the notation $\bfv \mapsto \g_\bfq(\bfv)$
to denote the isomorphism   from $T_\bfq\Q$ to $T^*_\bfq\Q$
induced by the scalar product $\g_{\bfq}[\cdot,\cdot]$.
Namely, for every $\bfv\in T_\bfq\Q$, the $1$-form
${\g}_{\bfq}(\bfv)$ is defined  by
\bel{gisdef}
\left\langle {\g}_{\bfq}(\bfv) , \w \right\rangle
~\doteq ~\g_{\bfq}[\bfv,\w]\qquad \forall \w\in T_\bfq\Q\,,
\eeq
where $\langle \cdot ,\cdot\rangle$ is the natural duality
between the tangent space $T_\bfq\Q$ and the cotangent
space $T^*_\bfq\Q$.


If $\bfq\in \Q$ and $W\subset T_{\bf q}$,~
$W^\perp$ denotes the subspace of $ T_{\bf q}$ consisting of all
vectors that are orthogonal to every vector in
$W$:
$$
W^\perp ~\doteq ~\left\{ \bfv\in T_{\bf q}\quad|
\quad \g_{\bf q}[\bfv, \w] = 0\quad\forall \w\in W\right\}.
$$

For a given distribution $E\subset T\Q$, the {\it orthogonal distribution} $E^\perp\subset T\Q$ is defined by setting $
E_{\bf q}^{\perp}\doteq (E_{\bf q})^\perp
$, for every $\bfq\in\Q$.


\item[{\bf 2)}]
Throughout the following,
we shall assume that the holonomic distribution $\Delta$ and
the non-holonomic distribution $\Gamma$ satisfy the  {\em tranversality condition}
\bel{transv}
\Delta_\bfq+\Gamma_\bfq ~= ~T_\bfq\Q\qquad\qquad \forall \bfq\in\Q\,.
\eeq
Notice that this is equivalent to
\bel{trans2}
\Delta^\perp\cap\Gamma^\perp~=~\{0\}\,,\eeq
and implies $\nu\leq N$.

      \item[{\bf 3)}] The mechanical
 system is subject to {\it forces}.
 In the Hamiltonian formalism, these are
      represented by  vertical vector fields on the
      cotangent bundle $T^*\Q$. We recall that, in a natural
      system
      of coordinates $(q,p)$, the fact that $\bf F$ is  {\it vertical}
       means that its $q$-component is zero, namely
       ${\bf F}=\sum_{i=1}^{N+M} F_i\frac{\partial}{\partial p_i}$.

\item[{\bf 4)}] The constraints   (\ref{nhol})
             and (\ref{hol}) are dynamically implemented
             by reaction forces  obeying

{\sc D'Alembert condition:}~
   {\it If $t\mapsto \bfq(t)$ is a trajectory which satisfies
   both  (\ref{nhol})
and (\ref{hol}), and $\bfR(t)$ is the constraint
reaction at a time $t$, then
\bel{dalembert} \bfR(t) ~\in ~\ker\Big(\Delta_{\bfq(t)}\cap
\Gamma_{\bfq(t)}\Big).\eeq}
In other words, regarding the reaction force $R(t)$ as an
element of the
cotangent space $T^*_{\bfq(t)}\Q$, one has
\bel{dal2}\langle \bfR(t)\,,~\bfv\rangle ~=~0\qquad\qquad\forall
\bfv~\in~\Delta_{\bfq(t)} \cap \Gamma_{\bfq(t)}~\subseteq ~
T_{\bfq(t)}\Q\,.
\eeq
\end{itemize}

\v
\subsection{Equations of motion}

For each $\bfq\in\Q$, we shall use $\g^{-1}_\bfq$
to denote the inverse of the isomorphism $\g_\bfq$ at (\ref{gisdef}).
Moreover we define the scalar product on the cotangent space
$\g^{-1}_\bfq[\cdot,\cdot]: ~T_\bfq^*\Q\times T_\bfq^*\Q~\mapsto~\R$
by setting
$$
\g_\bfq^{-1}[\p,\tilde\p]~ \doteq ~\g_\bfq\,[\g_\bfq^{-1}(\p), \,
\g_\bfq^{-1}
(\tilde \p)]\qquad \qquad\forall~~ (\p,\tilde\p)~\in ~
T_\bfq^*\Q\times T_\bfq^*\Q \,.
$$
For every $\bfq\in\Q$, we shall use $\H(\bfq,\cdot)$
to denote the Legendre transform
of the map $\bfv\to \T(\bfq,\bfv)$, so that
\bel{hamiltonian}
\H(\bfq,\p)~ =~ \frac{1}{2}\g^{-1}_\q[\p,\p]
\qquad\left( = \T(\bfq,{\g}^{-1}(\p))\right)
\qquad
\forall~ \p\in T_\bfq^{*}\Q\,.\eeq
The map $\H:T^{*}\Q\to\rr$  will be called the
{\it Hamiltonian corresponding to the kinetic energy $\T$}.

{}Let $(q)$ and $(u)$ be local coordinates on $\Q$ and $\U $, respectively, such that the domain of the chart $(q)$
is mapped by $\pi$ into the domain of the chart $(u)$.
Let $(q,p)$ the natural local coordinates  on $T^*\Q$
corresponding to the coordinates $(q)$. From
Nonholonomic Mechanics  it follows that, given a smooth control
$t\mapsto u(t)$ (here regarded as a time-dependent
holonomic constraint),  the corresponding motion
$t\mapsto(q,p)(t)$ on $T^*\Q$
verifies the relations
\bel{Ham}\left\{\begin{array}{rl}\dot q(t)
&=~\ds\frac{\partial H}{\partial p}(q(t),p(t))\, ,\cr &\cr
\dot p(t) &=~\ds
-\frac{\partial H}{\partial q}(q(t),p(t))
+ F(t,q(t),p(t)) + R(t)\,,
\cr &\cr
\pi\circ q(t) &=~u(t),\cr &\cr
p(t) &\in ~{g}_{q(t)}(\Gamma_{q(t)})\,,\cr &\cr
R(t)&\in ~\ker\Delta_{q(t)}+\ker\Gamma_{q(t)} \,.\end{array}\right.
\eeq
We have used
$H$, $F$, $R$, and $g$ to denote the local
expressions of $\H$, $\bf F$, $\bfR$, and $\g$, respectively.
For simplicity, we use the same notation for the distributions  $\Gamma,\Delta$
on
$\Q$ and their local expressions in coordinates.
Moreover, we use the notation $v\mapsto g_q(v)$
to denote the local expression of the isomorphism $\bfv\mapsto\g_\q(\bfv)$.

The first two equations   in (\ref{Ham}) are the dynamical equations
written  in Hamiltonian form. The first one yields the inverse of the Legendre transform.
The third equation represents the  (holonomic)
control-constraint. Relying on the fact that
$\pi$ is a submersion, we shall always choose local
coordinates $(q^r)$ and $(u^\alpha)$ such that $q^{N+\alpha} = u^\alpha$,
for all $\alpha=1,\dots,M$.   We thus regard this equation as prescribing
a priori the evolution of the last $M$ coordinates: $q^{N+1}, \ldots, q^{N+M}$.
The fourth relation   in (\ref{Ham}) is the Hamiltonian version
of the non holonomic constraint (\ref{nhol}).
The fifth relation  is clearly equivalent
to (\ref{dalembert}), i.e.~it  represents d'Alembert's condition.

\begin{remark} {\rm Although a global,  intrinsic
formulation of these equations can be given (see \cite{Marle})
we are here  mainly interested in their local coordinate-wise expression.
Indeed, a major goal of our analysis is to understand
the functional dependence on the time derivative $\dot u$ of the control.}
\end{remark}
The special case where
no active constraints are present can be obtained
by taking $\Delta\equiv T^*\Q$, i.e. $\ker(\Delta) =\{0\}$.  In this case,
(\ref{Ham}) reduces to the standard
Hamiltonian version of the dynamical equations
with non-holonomic constraints, namely
\bel{Hamno}\left\{\begin{array}{rl}\dot q(t)
&=~\ds\frac{\partial H}{\partial p}(q(t),p(t))\, ,\cr &\cr
\dot p(t) &=~\ds
-\frac{\partial H}{\partial q}(q(t),p(t))
+ F(t,q(t),p(t)) + R(t)\,,
\cr &\cr
p(t) &\in ~{g}_{q(t)}(\Gamma_{q(t)})\,,\cr &\cr
R(t)&\in ~\ker\Gamma_{q(t)} \,,\end{array}\right.
\eeq

\vs

\section{Orthogonal decompositions of the tangent and the cotangent bundles}
\setcounter{equation}{0}

To derive a set
of equations describing the constrained motion, it will be convenient
to decompose both the tangent
bundle $T\Q$ and the cotangent bundle $T^*\Q$
as direct sums of three suitable vector sub-bundles.
We recall that $\Q$ is a manifold of dimension $N+M$, while
$\Delta$ and  $\Gamma$ are distributions on $\Q$, having dimensions
$N$ and $N+M-\nu$, respectively. (In view of the transversality
condition (\ref{transv}) one has $\nu\leq N$.)

\subsection{Tangent bundle}
\begin{definition} For every $\bfq\in\Q$,  we define the following three subspaces of
$T_\bfq\Q$.
\bel{TJQ}
(T_\bfq\Q)_I~\doteq~\Delta_\bfq\cap\Gamma_\bfq\,,\qquad
(T_\bfq\Q)_{II}~\doteq~\Gamma_\bfq^\perp, \qquad
(T_\bfq\Q)_{III}~\doteq~ (\Delta_\bfq\cap \Gamma_\bfq)^\perp\cap \Gamma_\bfq\,.
\eeq
\end{definition}

\begin{figure}
\centering
\includegraphics[scale=0.40]{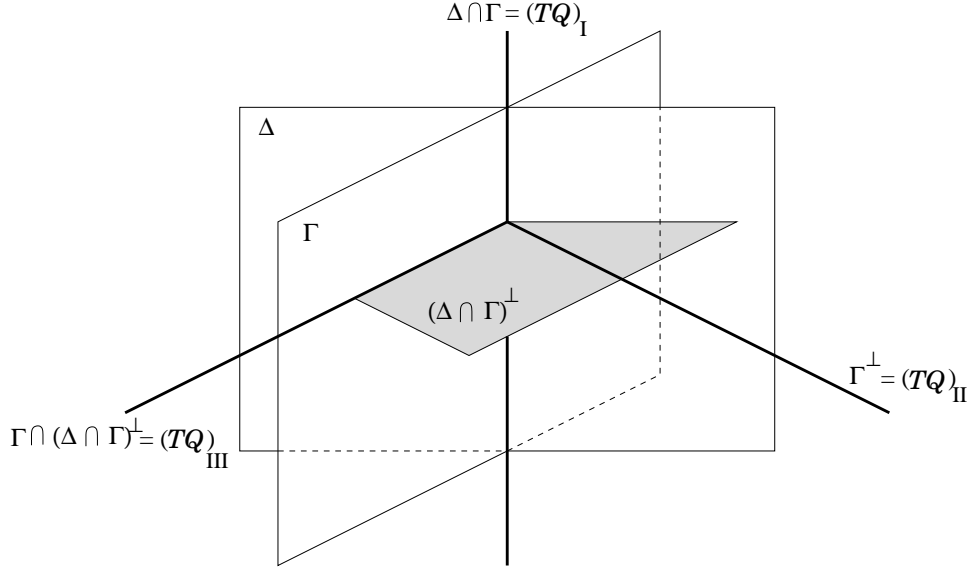}
\caption{\small{The orthogonal decomposition of $T\Q$.}}
\end{figure}

\begin{prop}\label{lemmasplit}
For each  $\bfq\in\Q$, the three subspaces in (\ref{TJQ}) are mutually orthogonal
and span the entire tangent space, namely
\bel{plustangent}
T_\bfq\Q ~=~ (T_\bfq\Q)_I\oplus (T_\bfq\Q)_{II}\oplus
(T_\bfq\Q)_{III}\,,
\eeq
Moreover,
 \bel{Gspan}
 (T_\bfq\Q)_{I}\oplus (T_\bfq\Q)_{III} ~=~ \Gamma\,.\eeq
If the transversality condition
(\ref{transv}) holds, then the above subspaces have dimensions
\bel{Tdim}
\hbox{\rm dim}(T_\bfq\Q)_{I}~ = ~N-\nu\,,
\qquad \hbox{\rm dim}(T_\bfq\Q)_{II}~ =~ \nu\,,
\qquad
\hbox{\rm dim}(T_\bfq\Q)_{III} ~=~ M\,.
\eeq
\end{prop}
\v
{\bf Proof.}  The orthogonality of the subspaces
$(T_\bfq\Q)_I$ and $(T_\bfq\Q)_{II}$ is immediately clear form the definitions.
Observing that
$$\Big((T_\bfq\Q)_{I}\oplus (T_\bfq\Q)_{II}\Big)^\perp ~=~
\Big((T_\bfq\Q)_{I}\Big)^\perp\cap \Big((T_\bfq\Q)_{II}\Big)^\perp~=~
 (\Delta_\bfq\cap \Gamma_\bfq)^\perp\cap \Gamma_\bfq~=~
 (T_\bfq\Q)_{III}\,,$$
 we obtain the orthogonal decomposition (\ref{plustangent}).
 In particular, this implies
$$ (T_\bfq\Q)_{I}\oplus (T_\bfq\Q)_{III} ~=~ \Big((T_\bfq\Q)_{II}\Big)^\perp~=~\Gamma_\bfq\,.$$

Finally, if (\ref{transv}) holds, then
$$\hbox{dim}(\Delta_\bfq\cap\Gamma_\bfq) ~=~\hbox{dim}(\Delta_\bfq) +
\hbox{dim}(\Gamma_\bfq) -\hbox{dim}(T_\bfq\Q)~=~
N + (N+M-\nu) - (N+M) ~=~N-\nu\,.$$
Moreover, dim$\Big((T_\bfq\Q)_{II}\Big)$ = dim$(T_\bfq\Q)$
 - dim$(\Gamma_\bfq)=
\nu$. The last equality in (\ref{Tdim}) now follows from (\ref{plustangent}).
\endproof
\v
For $J\in\{I,II,III\}$, the orthogonal projection
onto the subspace $(T\Q)_J$
will be denoted by
\bel{PJdef}\P_J: T_\bfq\Q~\mapsto ~(T_\bfq\Q)_{J}\,.\eeq

\begin{remark} {\rm Under the assumption (\ref{transv}),
the third subspace in the decomposition  (\ref{plustangent})
can be equivalently written as
$$(T_\bfq\Q)_{III}~=~\P^\Gamma(\Delta_\bfq^\perp),$$
where $\P^\Gamma$ denotes the orthogonal projection on the space
$\Gamma$.
Indeed, for any vector $\bfv\in T_\bfq\Q$ one has
$\bfv~\in~ \P^\Gamma(\Delta_\bfq^\perp) $  iff
$$\bfv~\in~ \Gamma_\bfq\cap (\Delta_\bfq^\perp + \Gamma_\bfq^\perp)
~=~ \Gamma_\bfq\cap  (\Delta_\bfq\cap \Gamma_\bfq)^\perp.
$$
}
\end{remark}

\subsection{Cotangent bundle}

Thanks to the isomorphism  ${\g}:T\Q\mapsto T^*\Q$ defined at
(\ref{gisdef}), one can use  (\ref{plustangent}) to obtain
a similar decomposition of the cotangent bundle as the
direct sum of three vector bundles:
\bel{pluscotangent}
T^*\Q ~= ~(T^*\Q)_I\oplus (T^*\Q)_{II}\oplus (T^*\Q)_{III},
\eeq
where, for  $J\in\{I,II,III\}$,
\bel{cotdec}
(T_\bfq^*\Q)_J~\doteq~{\g}_{\bfq}\Big((T_\bfq\Q)_J\Big)\,.
\eeq
We denote by
$\P^*_J: T^*_\bfq\Q\to (T^*_\bfq\Q)_{J}$ the orthogonal projection
w.r.t.~the metric $\g^{-1}$.
The above construction yields
\bel{P*}\P^*_J ~=~ \g\circ\P_J\circ\g^{-1}\,.\eeq
\v
\begin{remark} {\rm
By (\ref{P*}), the orthogonal
decompositions of the tangent and cotangent
bundles at (\ref{plustangent}) and (\ref{pluscotangent})
have the following property.
Recalling (\ref{gisdef}),
assume that $\p = \bfg( \dot \bfq)$, so that $\dot \bfq = \bfg^{-1}(\bfp)$.
Then,
 for $J\in \{ I, \, II, \,III\}$,
\bel{propro}\bfp_{\strut J}~\doteq~
\P^*_J(\bfp)~=~\bfg(\P_J(\dot \bfq))\,,
\qquad\qquad \dot\bfq_{\strut J}~\doteq
~\P_J(\dot \bfq) ~=~\bfg^{-1}(\P^*_J(\bfp))\,.\eeq
In other words, the $J$-th component of $\bfp$ depends only on the $J$-th component of $\dot \bfq$, and viceversa.}
\end{remark}

In view of Proposition \ref{lemmasplit}
one obtains:
\begin{prop}\label{colemmasplit} If the transversality condition
(\ref{transv})
holds, then
\bel{dim}
\hbox{\rm dim}\Big((T^*_\bfq\Q)_{I}\Big)~ = ~N-\nu\,,
\qquad \hbox{\rm dim}\Big((T^*_\bfq\Q)_{II}\Big)~ =~ \nu\,,
\qquad
\hbox{\rm dim}\Big((T^*_\bfq\Q)_{III} \Big)~=~ M\,.
\eeq
Moreover, for every $\q\in\Q$ the three subspaces
$(T^*_\bfq\Q)_{I}, (T^*_\bfq\Q)_{II}, (T^*_\bfq\Q)_{III}$
are pairwise orthogonal (w.r.t.~the metric $\g^{-1}$).
In particular,  $(T^*_\bfq\Q)_{II} = \ker(\Gamma_\bfq)$.
\end{prop}

\v
Next, we observe that by (\ref{transv}) the differential
of the submersion $\pi:\Q\mapsto U$ is one-to-one when restricted to
$(T_q\Q)_{III}$. More precisely,
$D\pi: (T_\bfq\Q)_{III}~\mapsto~T_{\pi(\bfq)}\U$
is a bijective linear map.  Its inverse will be denoted as
\bel{bfhdef}\bfh:T_{\pi(\bfq)}\U~\mapsto ~(T_\bfq\Q)_{III}\,.\eeq
For any
$ \bfv \in T_{\pi(\bfq)}\U$, by $\bfh( \bfv)$ we thus
denote the  unique vector in $(T_\bfq\Q)_{III}$
such that
\bel{hudef}
D\pi\cdot (\bfh(\bfv)) ~=~\bfv\,.\eeq
Moreover,  we write
\bel{bfkdef}
\bfk(\bfv) ~\doteq~\bfg_\bfq(\bfh(\bfv))~\in~(T_\bfq^*\Q)_{III}\,.\eeq
  We claim that
the vector $\bfh(\bfv)$ defined above can be characterized as
the unique vector $\bfw\in \Gamma_\bfq$ where the following constrained minimum is attained:
\bel{argmin}\bfh(\bfv)~=~\hbox{arg}\!\!\!\!\!\!\!\!\!\!\!\!\min_{\strut
\bfz\in \Gamma_\bfq\,,~D\pi
\cdot\bfz= \bfv}~ \bfg_\bfq[\bfz,\bfz]\,.\eeq
Indeed, since the Riemann metric $\bfg_\bfq$ is positive definite, the
right hand side of (\ref{argmin}) is well defined.
Consider any vector $\bfw\in T_\bfq\Q$ such that
\bel{wmin} \bfw\in \Gamma_\bfq\,,\qquad
 D\pi\cdot \bfw ~=~\bfv\,.\eeq
The above equalities imply
$$\bfw~=~\bfh(\bfv)+\bfw_{I}$$
for some $\bfw_{I}\in (T_\bfq\Q)_I$.
Since the vectors  $\bfh(\bfv)\in (T_\bfq\Q)_{III}$
and $ \bfw_{I}$
are orthogonal, we have
$$\bfg_\bfq[\bfw,\bfw]~=~\bfg_\bfq[\bfh(\bfv),\, \bfh(\bfv)]+\bfg_\bfq[\bfw_I,\bfw_I]~
\geq~\bfg_\bfq[\bfh(\bfv),\, \bfh(\bfv)].$$
This proves our claim.

\begin{remark}\label{r3.3} {\rm  Let $\bfu\in\U$ and choose $\bfq\in\Q$ such that
$\bfu=\pi(\bfq)$.
Identifying the tangent spaces $T_\bfu\U\approx T_\bfq\Q/\Delta_\bfq$,
the inner product
$$\langle \bfv,\bfv'\rangle_{\bfu}~\doteq~\bfg_\bfq[ \bfh(\bfv),\,\bfh(\bfv')] $$
can be seen as a Riemann metric on the space of leaves.   In general,
this metric is not canonically defined, because
it depends on the choice of a particular point $\bfq$ on each leaf.
{}From a more precise result proved in
Section~7 it will follow that, if this metric is independent of the choice of
$\bfq$,
then the system is ``fit for jumps". Namely, the right hand sides of the
dynamical equations (\ref{c1}) are
affine functions of the time derivative $\dot \bfu$ of the control.  As a consequence,
trajectories of the system can be meaningfully defined also in connection with discontinuous
control functions \cite{B-R1, B-R2}.
}
\end{remark}
\v
\section{A closed system of control equations}
\setcounter{equation}{0}

In the original formulation (\ref{Ham}), the motion
is characterized in terms of a family of ODEs
coupled with a set of constraints.
Relying on the decompositions at
(\ref{plustangent}) and (\ref{pluscotangent}),
we will show that the motion can be described by
a system of differential equations.

The following result, providing different ways
to express the
constraint (\ref{nhol}), is straightforward.
\v
 \begin{lma}\label{nholplemma}
Let  $t\mapsto \q(t)\in \Q$ be a $\C^1$ map.
Using the same notation as in
(\ref{propro}), the non-holonomic constraint $\dot \bfq(t)\in
\Gamma_{\bfq(t)}$ can be expressed in
any of the following equivalent forms:
\bel{equiv}
\p(t)\in \g(\Gamma_{\bfq(t)})~~~
\Longleftrightarrow ~~~\dot \bfq_{II}(t) = 0
 ~~~\Longleftrightarrow~~~ \p_{II}(t) = 0.\eeq
\end{lma}
\v
On  a given, natural chart $(q,p)$ on $T^*\Q$, and for
$J\in \{I,\, II,\, III\}$,
we denote by $P_J^*$  the ($q$-dependent) matrix
representing the projection
$\P_J^*$.
Let $g = (\bfg_{r,s})$ be the matrix representing
the Riemannian metric $g$
 in the $q$-coordinates. In turn, the inverse matrix
 $g^{-1} = (g^{r,s})$ represents the metric $g^{-1}$
 on the cotangent space. In this coordinate system,
the components of the force and of the constraint reaction
will be denoted as
$F_{J} \doteq P^*_J \,F$ and $R_{J} \doteq
P^*_J \,R$.   The
coordinate representations of the linear maps $\bfh$ and
$\bfk$ introduced at (\ref{bfhdef})-(\ref{bfkdef}) will be denoted by
$h$ and $k$, respectively.

Throughout the following, we use a local system of
{\it adapted} coordinates, so that $q^{N+\alpha}= u^\alpha$
for $\alpha = 1,\ldots,M$.   If the holonomic active constraint
is satisfied,  by (\ref{propro}) and (\ref{hudef})-(\ref{bfkdef})
this
implies the identity
\bel{dqu}P_{III}(\dot q(t))~ = ~h(\dot u(t))\,,\qquad\qquad
P^*_{III} (p(t))~=~k(\dot u(t))\,.\eeq

In terms of the orthogonal decompositions (\ref{plustangent})
and (\ref{pluscotangent}), the non-holonomic and the active
holonomic constraints yield  the system of equations
\bel{c23}\left\{\begin{array}{l}
\dot q_{II}(t) =~0  \\\,\\
\dot q_{III}(t) =~h(\dot u(t))\end{array}\right.
 \qquad\qquad\left( \hbox{which holds if and only if}
 \quad
\left\{ \begin{array}{l}
p_{II}(t) =~0 \\\,\\
 p_{III}(t) =~k(\dot u(t))\end{array}\right.\right).
\eeq

To complete the description of the motion, it remains to
derive the equations for $\dot q_I$ and $\dot p_I$.
By (\ref{propro}) it follows
\bel{dqeq} \dot q_I ~=~P_I(\dot q) ~=~P_I ( g^{-1}(p))~=~
g^{-1}(P^*_I(p))\,.\eeq
Moreover, differentiating  $p_I = P_I^*(p)$ w.r.t.~time
and using (\ref{Ham}), we obtain
\bel{pIdot}\bega{rl}
\dot p_I&\ds =~\left({\partial P_I^*\over\partial q}
\cdot \dot q\right) (p)
+ P_I^*(\dot p)  \cr&\cr
&\ds=~\left({\partial P_I^*\over\partial q}
\cdot \Big(g^{-1}(p_I) +h(\dot u)\Big)\right) (p_I + p_{III}) +
P_I^*\left(-\frac{\partial H}{\partial q}(q(t), p(t)) + F(t,q(t),p(t)) + R(t)\right)\cr &\cr
&\ds=~\left({\partial P_I^*\over\partial q}
\cdot g^{-1}\Big(p_I +k(\dot u)\Big)\right) (p_I + k(\dot u))
-\frac{1}{2}
P_I^*\left(\frac{\partial g^{-1}}{\partial q}[p_I+k(\dot u)\,,~ p_I+k(\dot u)]\right)+F_I\,.
\enda
\eeq
Indeed, the fifth relation in (\ref{Ham}), i.e.~
d'Alembert's condition, implies $P^*_I(R(t)) = 0$.

The first two terms on the right hand side of (\ref{pIdot})
can be written in a simpler form, introducing the bilinear
($q$-dependent, possibly not symmetric), $\R^{N+M}$-valued
map
\bel{thIdef}
(p,\tilde p)
~\mapsto~ \theta_I[p,\tilde p] ~\doteq ~\left({\partial P_I^*\over\partial q}
\cdot g^{-1}(p)\right) (\tilde p) -\frac{1}{2}
P_I^*\left(\frac{\partial g^{-1}}{\partial q}[p,\, \tilde p]\right).
 \eeq
Setting
\bel{centrifugal}
 \Upsilon[p_I\,,\,\dot u]~\doteq
 ~\theta_I[p_I\,,k(\dot u)] + \theta_I[k(\dot u)\,,~p_I],\qquad
 \qquad\Psi[\dot u,\,\dot u]~\doteq
 ~\theta_I[k(\dot u)\,,k(\dot u)],
\eeq
the equations
(\ref{dqeq}) and (\ref{pIdot}) can be written in the form
\footnote{We recall  that $\dot q_{I}(t)\equiv P_I(\frac{dq}{dt})$ {\it while} $\dot p_{I}(t) \equiv \frac{d}{dt}(P^*_I(p))$}
\bel{c1}
\left\{ \bega{rl}
\dot q_{I}(t) &=~g^{-1}(p_I)\,,\cr &\cr
\dot p_{I}(t)&=~\theta_I[p_I, p_I] + \Upsilon[p_I, \dot u] +
\Psi[\dot u, \dot u] + F_I\,.\enda\right.
\eeq
We can now state the main result of this section.

\begin{thm}  \label{t:42}
 Let $\I$ be a time interval and let
$u:\I\to\R^M$
be a $\C^1$ control function. Let $q:\I \mapsto \R^N$ be a
$\C^1$ path and set $p(t)\doteq g(\dot q(t))$.
Then the path $(q,p):\I\to T^*\Q$ satisfies the non holonomic
equations (\ref{Ham}) (with constraints as controls)
if and only if  the following two conditions are satisfied.
\begi
\item[(i)] For $J\in \{I,II,III\}$,
the components $\dot q_J \doteq P_J (\dot q(t))$ and
$p_J(t)\doteq P^*_J(p(t))$ satisfy the equations
in (\ref{c23}) and (\ref{c1}).
\item[(ii)] At some time $t_0\in \I$
one has
$q^{N+\alpha}(t_0) = u^\alpha(t_0)$ for all $\alpha=1,\dots, M$.
\endi
 \end{thm}
\v

 {\bf Proof.} {\bf 1.}  By the previous analysis, if
 all relations in (\ref{Ham}) are satisfied, then the
 equations (\ref{c23}) and (\ref{c1}) hold.  Moreover,
 the condition $\pi\circ q(t) = u(t)$ implies (ii), for all
 times $t\in \I$.
 \v
{\bf 2.} Next, assume that all the equations in (\ref{c23}) and (\ref{c1}) hold, and $\pi(q(t_0)) = u(t_0)$ at some time $t_0\in \I$.
By  (\ref{equiv}), the equation $\dot q_{II}=0$
implies the fourth relation in (\ref{Ham}).

For any $t\in \I$, if $\pi(q(t)) - u(t)=0$ then
$$\frac{d}{dt}\Big( \pi(q(t))-u(t)\Big) ~
=~D\pi\cdot \dot q(t)- \dot u(t) ~=~
D\pi\cdot \dot q_{III}(t) - \dot u(t)~=~0\,.$$
{}From
the assumption $\pi(q(t_0)) - u(t_0)=0$ and
the regularity of the coefficients of the equation, we conclude that
$\pi(q(t)) - u(t)=0$ for all $t\in\I$.  Hence the third relation
in (\ref{Ham}) holds.

The equations for $\dot q$ in (\ref{c23}) and (\ref{c1}) imply
$$\dot q ~=~\dot q_I + \dot q_{II} + \dot q_{III} ~=~g^{-1}(p_I) +h(\dot u)~=~
g^{-1}(p_I +p_{II}+ p_{III})~=~\frac{\partial H}{\partial p}
(q, p)\,,$$
This yields the first equation in (\ref{Ham}).

Finally, by defining
\bel{react2}
R(t)~\doteq ~
\dot p + \frac {\partial H}{\partial q}(q,p)
- F(t,q(t),p(t)),\eeq
the second equation in (\ref{Ham}) is clearly satisfied.
It remains to check  the fifth equation in (\ref{Ham}), namely
$R(t)\in \ker\Delta_q(t) +  \ker\Gamma_q(t)$. Since  by construction
$\ker\Delta +  \ker\Gamma = (T^*\Q)_{II} + (T^*\Q)_{III}$,
 this is  equivalent to proving that  $R_I = 0$.
Using the second equation in (\ref{c1}) we find
$$\bega{rl}
R_I &\ds
=~P_I^*(R) ~=~P_I^*(\dot p) + P_I^*
\left( \frac {\partial H}{\partial q}(q,p)\right)
- F_I\cr &\cr
&\ds =~\dot p_I - \left(\frac{\partial P^*_I}{\partial q}\cdot \dot q\right)(p)+ P_I^*
\left( \frac {\partial H}{\partial q}(q,p)\right)
- F_I\cr
&\cr
&=\ds~ \dot p_I -\left({\partial P_I^*\over\partial q}
\cdot g^{-1}\Big(p_I +k(\dot u)\Big)\right) (p_I + k(\dot u)) +
\frac{1}{2}
P_I^*\left(\frac{\partial g^{-1}}{\partial q}[p_I+k(\dot u)\,,~ p_I+k(\dot u)]\right)-F_I\cr &\cr
&=~0\,.
 \enda
 $$
 \endproof

\section{The equations in $\Delta$-adapted coordinates}\label{reducedsect}
\setcounter{equation}{0}

Let $q=(q^1,\dots,q^{N+M})$ be $\Delta$-adapted coordinates
on a open set $\mathcal{O}\subseteq\Q\,$, so that
$$\Delta~=~span \Big\{\frac{\partial}{\partial q^1}\,,\, \dots\,,\,
\frac{\partial}{\partial q^N}\Big\}.$$
The main goal of this section is to write  the equations of motion
explicitly in terms of these state  coordinates, together
with adjoint coordinates
$\xi_j$, $1\leq j\leq N+M$ corresponding
to suitable  bases
 of the cotangent
 bundle $T^*\Q$, decomposed as in (\ref{pluscotangent}).
 It will turn out that  the relevant equations will involve only the $2N-\nu$ variables $q^1,\dots,q^N,\xi_1,\dots,\xi_{N-\nu}$.
\v

Consider a family $\{\V_1,\,\dots\,,\V_{N+M}\}$ of smooth,
linearly independent
vector fields on  $\Q$, such that, at least on the open set $\O$,
\bel{familyVF}\begin{array}{l}
(T_\bfq\Q)_I~=~span\Big\{\V_1(\bfq),\ldots,\V_{N-\nu}(\bfq)\Big\}\,,
\\\,\\
 (T_\bfq\Q)_{II} ~=~ span\Big\{\V_{N-\nu+1}(\bfq),\ldots,\V_{N}(\bfq)
 \Big\}\,,
 \\\,\\
 (T_\bfq\Q)_{III}~ =~ span\Big\{\V_{N+1}(\bfq),\ldots,
 \V_{N+M}(\bfq)\Big\}.
 \end{array}
\eeq
Throughout the following, we
assume that the  vectors $\{\V_1,\,\dots\,,\V_{N-\nu}\}$, which generate $(T_\bfq\Q)_I$, are mutually orthogonal, i.e.
\bel{orthI}
\g[\V_r,\V_s] ~=~ 0
\qquad\forall r,s\in \{1,\dots,N-\nu\},\quad r\not= s\,.
\eeq
In addition, for
$i=1,\dots, N+M$, we define the basis of  cotangent vectors
\bel{Forme}
{\bf \Omega}_i~\doteq ~\ds \frac{ \g(\V_i)}{\g[\V_i,\,\V_i]},
\eeq
By (\ref{cotdec}), this yields
$$\begin{array}{l}
(T^*_\bfq\Q)_I~=~span\Big\{{\bf \Omega}_1(\bfq),\dots,
{\bf \Omega}_{N-\nu}(\bfq)\Big\}\,,\\\,\\
 (T^*_\bfq\Q)_{II} ~=~  span\Big\{{\bf \Omega}_{N-\nu+1}(\bfq),
 \dots,{\bf \Omega}_N(\bfq)\Big\}\,,
 \\\,\\(T^*_\bfq\Q)_{III} ~=~  span\Big\{{\bf \Omega}_{N+1}(\bfq),
 \dots,{\bf \Omega}_{N+M}(\bfq)\Big\}\,.
 \end{array}
$$
By (\ref{orthI}) and (\ref{Forme}), the differential forms
$\{{\bf \Omega}_1,\,\dots\,,{\bf \Omega}_{N-\nu}\}$ are mutually
orthogonal with respect to the metric $\g^{-1}$,
namely
\bel{orthII}
\g^{-1}[{\bf \Omega}_r,\,{\bf \Omega}_s] ~=~ 0
\qquad\forall r,s\in \{1,\dots,N-\nu\},\quad r\not= s\,.
\eeq
Moreover,  the basis
$\Big\{{\bf \Omega}_1(\bfq),\dots,{\bf \Omega}_{N-\nu}(\bfq)\Big\}$
is dual to the basis
$\Big\{\V_1(\bfq),\ldots,\V_{N-\nu}(\bfq)\Big\}$, i.e.
$$
{\bf \Omega}_r(\bfq)\Big(\V_s(\bfq)\Big) ~= ~\delta_{r,s}  \qquad\forall r,s=1,\dots,N-\nu.
$$
This choice of orthogonal bases makes it easy to compute
the projections $\P_I$ and $\P_I^*$.
Indeed, for any tangent vector $\bfw$ and any cotangent
vector $\bfp$ one has
\bel{projections}
  \P_I (\bfw)~=
   ~\sum_{\ell=1}^{N-\nu}\ds\frac{\g[\bfw, \V_\ell]}{\g[\V_\ell,\V_\ell]} \V_\ell,
  \eeq
\bel{coprojections}
  \P^*_I (\p)~=
  \sum_{\ell=1}^{N-\nu}
 {\g^{-1}[\p, {\bf \Omega}_\ell]\over \g^{-1}[{\bf \Omega}_\ell,\,
 {\bf \Omega}_\ell]}\,
  {\bf \Omega}_\ell\,.   \eeq

In the following, to simplify notation, whenever repeated indices taking values
from 1 to $N+M$ are summed,
the summation symbol will be omitted.
On the other hand, summations ranging over a smaller set of indices
will be explicitly written. As before, we
let $g=(\bfg_{r,s})$ be the matrix representing the
Riemannian metric  $\g$ in the $q$-coordinates. In turn,
the inverse matrix
$g^{-1}=(g^{r,s})$  represents the metric $\g^{-1}$
on the cotangent space.

 For $\ell=1,\dots,N+M$, let
   ${V_\ell}^1,\ldots,{V_\ell}^{N+M}$ be the
   $q$-components of
 $\V_\ell$, so that $\V_\ell={V_\ell}^i
 \frac{\partial}{\partial q_i}$.
If the  (column) vector  $w= (w^s)\in\R^{N+M}$ yields
the coordinate representation of $\bfw$, then
 by (\ref{projections}) the projected vector  $\P_{I}(\bfw)$
is a column vector with coordinates
 $
  (P_{I}\,w)^s~=~(P_I)^s_r w^r$,
  where
  the $(N+M)\times (N+M)$ matrix $P_{I}$ is defined by
 $$
 (P_{I})^s_r ~\doteq ~\sum_{\ell=1}^{N-\nu}\frac{\bfg_{r,k}{V_\ell}^k
 {V_\ell}^s}{\bfg_{i,j}{V_\ell}^i
 {V_\ell}^j}\qquad \qquad r,s=1,\dots,N+M\,.
  $$

Similarly,
let $\Omega_{r,1},\dots,\Omega_{r,N+M}$
 be the components of  $\Om_r$, so that
 ${\bf \Omega}_r=\Omega_{r,s}dq^s$.
If the (row) vector   $p\in\R^{M+N}$
yields
 the coordinate representation
 of the covector $\p$, then by  (\ref{coprojections})
 the projected the vector  $\P^*_{I}(\p)$ has coordinates given by
$(p\,P^*_{I})_s = p_r (P^*_I)^r_s$,
  where
  the $(N+M)\times (N+M)$ matrix $P^*_{I}$ is defined by
\bel{copromat}
(P_{I}^*)^r_s~ \doteq~\sum_{\ell=1}^{N-\nu}\frac{g^{r,k}\,\Omega_{\ell,k}\,
 \Omega_{\ell,s}}{g^{i,j}\,\Omega_{\ell,i}\Omega_{\ell,j}}\,.
 \eeq

In order to write the second set of
equations in (\ref{c23}),
we need an explicit expression
of the ($q$-dependent) matrices $h$ and $k$.
Let us define the $(N+M)\times M$ matrix $V_{III}$ and the $M\times M$ matrix ${V_{III}}^{III}$ by setting
$$
V_{III} \doteq \Big({V_{N+\alpha}}^r\Big),\qquad {V_{III}}^{III} \doteq \Big({V_{N+\alpha}}^{N+\beta}\Big)
$$
Here and in the sequel,  Greek
 indices such as $\alpha,\beta$ range from $1$ to $M$, while Latin indices such as $r,s$
 range from 1 to $N+M$.
Notice that, by (\ref{familyVF}), the columns of $V_{III}$ span
$(T_q\Q)_{III}$.

 Recalling the definitions (\ref{hudef})-(\ref{bfkdef}),
 and the identities $q^{N+\alpha} = u^\alpha$,
it is easy to check
that the injective linear map ${\bf h}:T\U_{\pi(\bfq)}\mapsto
(T_\bfq\Q)_{III}\subset T_\bfq\Q$
is represented by the $(N+M)\times M$ matrix $$h~ =~V_{III}\cdot  \Big({V_{III}}^{III}\Big)^{-1}.$$
In turn, the linear map
${\bf k}(\cdot)~ = ~\g({\bf h}(\cdot))$ is represented by the
$(N+M)\times M$ matrix $$k\doteq g\cdot h ~=~
  g\cdot V_{III}\cdot  \Big({V_{III}}^{III}\Big)^{-1}.$$
  \v
Using this particular system of coordinates,  the equations
of motion (\ref{c23}), (\ref{c1})
for the non-holonomic system with active constraints
 can be written in the following form.

\begin{prop} Let   $(q,p, u)(\cdot)$ be as in Theorem \ref{t:42}.
Moreover, let $(\xi,\, \eta,\, \lambda)$ be the components
of $\p=p_i dq^i$ w.r.t.~the frame
$\{{\bf \Omega}_1,\dots,{\bf \Omega}_{N+M}\}$, so that
   $$
   \p(t)~ = ~\sum_{\ell=1}^{N-\nu}\xi_\ell(t){\bf \Omega}_\ell +  \sum_{\ell=N-\nu+1}^{N}\eta_\ell(t){\bf \Omega}_\ell
   + \sum_{\ell=N+1}^{N+M}\lambda_\ell(t){\bf \Omega}_\ell\,.
   $$
   Then the curve  $t\mapsto
   (q(t),p(t),u(t))$ provides a solution to the system
  (\ref{Ham}) if and only if
   the curve  $t\mapsto (q(t),\xi(t), \eta(t), \lambda(t), u(t))$
   is a solution of
   \bel{Hamclosedx1}\left\{\begin{array}{rll}
 \dot q^r &=~\ds
 \sum_{\ell=1}^{N-\nu}\frac{{V_\ell}^r}{\bfg_{i,j}{V_\ell}^i{V_\ell}^j}
 \,\xi_\ell + \sum_{\alpha=1}^M
 {h_{\alpha}}^r\, \dot u^{\alpha} \qquad &r=1,\dots,N\,,
\\\,\\
q^{N+\alpha} &= ~  u^\alpha  \qquad& \alpha=1,\dots,M \,,
\cr&\cr
 \dot \xi_m &=~\ds \tilde\theta_m [\xi,\, \xi] +
 \Tilde\Upsilon_m\,[\xi\,,\dot u] +
  \Tilde \Psi_{m}
 [\dot u ,\dot u] +  \Tilde F_m\qquad &m=1,\dots,N-\nu\,,\\\,\\
 \eta&=~\ds 0\,, &\\\,\\
 \lambda &= ~\ds(V_{III})^T\cdot k\cdot\dot u \,, &
 \end{array}
 \right. \eeq
where the superscript $^T$ denotes transposition and,
recalling (\ref{thIdef})-(\ref{centrifugal}),
for every $m=1,\dots,N-\nu$ we set
\bel{thetatilde}
\begin{array}{c}
\tilde\theta_m [\xi,\hat \xi]\doteq \sum_{a,b=1}^{N-\nu}\tilde\theta^{a,b}_m\,\xi_a\hat \xi_b,\quad\hbox{with} \quad \tilde\theta^{a,b}_m ~\doteq~\left( [{\theta_I}]_i^{r,s}\Omega_{r,a}\Omega_{s,b} -\ds\sum_{j=1}^N \frac{\partial\Omega_{i,a}}{\partial q^j}g^{j,s}\Omega_{s,b}\right){V_m}^i,\\\,\\

\Tilde\Upsilon_m\,[\xi\,,\dot u]~\doteq~
\sum_{a=1}^{N-\nu}\Tilde\Upsilon_{\alpha,m}^a\,\xi_a\,\dot u^{\alpha},
\quad\hbox{with} \quad\Tilde \Upsilon_{\alpha,m}^a~ \doteq~
\left(\Upsilon_{\alpha,i}^r\,\Omega_{r,a}
-\ds  \frac{\partial\Omega_{i,a}}{\partial q^{N+\alpha}}
 - \ds \sum_{j=1}^N \frac{\partial\Omega_{i,a}}
 {\partial q^j} {h_{\alpha}}^j\right){V_m}^i,\\\,\\

  \Tilde \Psi_{m} \,
 [w, \tilde w]~\doteq ~\Tilde \Psi_{\alpha,\beta,m} \,
 w^{\alpha} \tilde w^{\beta},
 \quad\hbox{with} \quad \ds\Tilde\Psi_{\alpha,\beta,m}~\doteq~ \Psi_{\alpha,\beta,i}{V_m}^i\,,\qquad\qquad
 \Tilde F_m ~\doteq~  F_i {V_m}^i\,.
 \end{array}
\eeq
\end{prop}
Indeed, (\ref{thetatilde}) provides the representation
of the maps $\theta_I, \Upsilon$, and $\Psi$ in
(\ref{thIdef})-(\ref{centrifugal}), with the present choice
of coordinates.

\v
\section{Systems which are ``fit for jumps"}
\label{s:4}
 \setcounter{equation}{0}
We now examine the connections between the following
properties:
\begin{itemize}
\item[i)] The continuity of the control-to-trajectory
map $u(\cdot)\mapsto (q(\cdot),p(\cdot))$.
\item[ii)] The vanishing of the ``centrifugal" term $\Psi[\dot u,\dot u]$ in (\ref{centrifugal}).
\item[iii)]
The invariance of the distribution
$\Gamma\cap(\Gamma\cap\Delta)^\perp$
in $\Gamma$-constrained inertial motions.
\end{itemize}

We recall that, if $\Psi[\dot u, \dot u]\equiv 0$, then the equations
of motion (\ref{Hamclosedx1}) are affine w.r.t.~the time derivative
$\dot u$ of the control function. In this case, the  control-to-trajectory
map $u(\cdot)\mapsto (q(\cdot),p(\cdot))$ can be continuously  extended to a family of
possibly discontinuous control functions $u(\cdot)$ with bounded variation \cite{B-R1}.
As in \cite{AB1},
we then say that the system is {\em fit for jumps.}  On the other hand, if the quadratic term
 $\Psi[\dot u,\dot u]$ does not vanish, this means that small vibrations of the control
produce a centrifugal force. As shown in \cite{B-R4}, this feature can be exploited
as an additional way to control or stabilize the system.
\v
As in Section~2,
$\T$ will denote the kinetic energy associated with
the metric $\g$, defined at (\ref{kinetic}), while $\Gamma$
denotes a non-holonomic distribution on the manifold $\Q$.

\begin{definition} Let $\I$ be a time interval.
A $\C^2$ map $\q:\I\to\Q$ will be called a
{\em free} {\em inertial motion} if, in any set of coordinates, its local expression $q(\cdot)$ provides a solution to the Euler-Lagrange equations
 \bel{geo}
     \frac{d}{dt}\frac{\partial {\mathcal T}}
     {\partial \dot q} -
     \frac{\partial {\mathcal T}}{\partial q} ~=~ 0.
     \eeq

     \end{definition}

     \begin{definition}\label{im} A $\C^2$ map $q:\I\to\Q$
     will be called a {\em $\Gamma$-constrained
     inertial motion}  if, in any set of coordinates, its local expression $q(\cdot)$  is a solution of
 \bel{constraintgeo}
      \frac{d}{dt}\frac{\partial {\mathcal T}}{\partial \dot q}
      - \frac{\partial {\mathcal T}}{\partial q} ~\in ~\ker\,\Gamma
     \eeq
     and satisfies the non holonomic constraints
     \bel{nhvel}
     \dot q ~\in~ \Gamma.
     \eeq
     \end{definition}

      \begin{definition}\label{freeinvariant}
      Let $S\subseteq T\Q$ be any set. We say that
      $S$ is {\em invariant for the  free inertial flow}, {\rm (}or
     equivalently: {\em inertially invariant)},
     if, for every free inertial motion  $\q:\I\to\Q$ such that
     $(\q,\dot \q)
      (t_0)\in S$ for some $t_0\in \I$, one has
      \bel{stareins0}(\q,\dot \q)(t)\in S\qquad \forall t\in \I.\eeq
     \end{definition}

     \begin{definition}\label{invariant}
     Let $S\subseteq\T\Q$ be any set. We say that $S$
     is {\em invariant for the  $\Gamma$-contrained
     inertial flow}, {\rm (}or
     equivalently: {\em $\Gamma$-inertially invariant)},
       if, for every $\Gamma$-constrained inertial motion (see Def. \ref{im})
   $\q:\I\to\Q$ such that  $(\q,\dot \q)
      (t_0)\in S$ for some $t_0\in \I$, one has \bel{stareins}
        (\q,\dot \q)(t)\in S\qquad \forall t\in \I.\eeq
     \end{definition}

In the next theorem we state a characterization for the $\Gamma$-inertially invariance of set $(T\Q)_{III}
    =\Gamma\cap(\Gamma\cap\Delta)^\perp$.
\v
\begin{thm}\label{udotlinear} The following conditions are equivalent:
\begin{itemize}
\item[\rm (1)] If  $F\equiv 0$, for every $\C^1$ control
$u(\cdot)$  and every solution $q(\cdot)$ of
the control differential equation (on $\Q$)
    $$
    \dot q~ = ~{ h} (\dot u),
    $$
    the map $$t~\mapsto~(q(t),p(t))~\doteq~ \Big(q(t)\,,\,
    {k}(\dot u(t))\Big)$$ is a solution of {\rm (\ref{Ham})}.
    In particular, $p_I(t)=0$ and hence $\dot q_I(t)=0$
    for all $t\in \I$.
    \item[\rm (2)] For every local chart $(q)$ the (vector-valued)  quadratic form $\Psi$ defined at (\ref{centrifugal})
    vanishes identically.
    \item[\rm (3)] For every local chart $(q)$, the (vector-valued) quadratic
    form $ v\mapsto\theta_I[v,v]$ in
    {\rm(\ref{thIdef}) } vanishes on the
    subspace $P_{III}^*(\R^{N+M})$.

    \item[\rm (4)] The sub-bundle $(T\Q)_{III}
    =\Gamma\cap(\Gamma\cap\Delta)^\perp$ is
    $\Gamma$-inertially invariant.
        \end{itemize}
        \end{thm}

{\it Proof.}
We first prove the implication (1)$\Longrightarrow$(4).
Let (1) hold, and let $q:\I\mapsto\Q$ be a
$\Gamma$-constrained inertial motion such that $\dot q_I(t_0)=0$.
Choosing $u(t) = \pi (q(t))$, it is clear that the equations
(\ref{constraintgeo})-(\ref{nhvel}) imply (\ref{Ham}),
with $F\equiv 0$.  Hence, if (1) holds, then $\dot q_I(t)\equiv 0$
for all $t\in\I$.   By construction,
the non-holonomic constraint yields $\dot q_{II}(t)\equiv 0$.
Therefore
$\dot q(t) = \dot q_{III}(t)\in (T_q\Q)_{III}$ for all $t\in\I$.
\v
Next, we prove the implication (4)$\Longrightarrow$(2).
Fix any point $q(0)=q_0$ and any vector $\omega_0\in \R^M$.
Let $t\mapsto q(t)$ be the $\Gamma$-constrained inertial
motion starting at $q_0$, with $\dot q(0) = h(\omega_0)$.
Then $q(\cdot)$ provides a solution to the system
(\ref{Ham}) in connection with the control
$u(t) \doteq \pi(q(t))$.
If (4) holds, then $\dot q_I(t)\equiv 0$ for all $t$, and hence
$p_I(t)\equiv 0$ as well.
On the other hand, our construction implies
$\dot u(0) =\omega_0$.    Therefore,
at time $t=0$, the second equation in (\ref{c1}) yields
$$0=\dot p_I(0) ~=~\Psi[\dot u(0),
\dot u(0)] ~=~\Psi[\omega_0,\omega_0]\,.$$
Since $q_0$ and $\omega_0$ can be chosen arbitrarily, we conclude that
$\Psi\equiv 0$, proving (2).
\v
We now prove the implication (2)$\Longrightarrow$(1).
Let a continuously differentiable control $u(\cdot)$ be given.
If condition (2) holds true and $F\equiv 0$, then
in the ODE (\ref{c1}) for $\dot p_I$ each term on the right hand side
is a linear or quadratic w.r.t.~$p_I$.
Therefore, if $p_I(t_0) =0$
at some time $t_0$ then $p_I(t)\equiv 0$
for all times $t\in I$, proving (1).
\v
Finally, we show that (2)$\Longleftrightarrow$(3).
Since the linear map $\bfk:T_{\pi(\bfq)}\U\mapsto
(T_\bfq\Q)_{III}$ is a bijection,
the equivalence of (2) and (3) follows directly by the
definition of $\theta_I$ and $\Psi$.
\endproof
\v



\begin{remark}{\rm
In the holonomic case
the the quadratic term in $\dot u$, which prevents the system to be ``fit for jumps",  is completely determined
by the relation between the kinetic metric and the
distribution  $\Delta$. In particular this quadratic term
accounts for the curvature of the
orthogonal distribution $\Delta^\perp$
(see \cite{B-R4} for a precise result in this direction).
On the other hand, when an additional non holonomic constraint
$\dot q\in \Gamma$ is present, the orthogonal distribution $\Delta^\perp$
can have zero curvature\footnote{We recall that $\Delta^\perp$
has zero curvature if the following holds.
Let $s\mapsto q(s)$ be any geodesic of the manifold
$\Q$ which intersects perpendicularly one of the leaves of the foliation
generated by the integrable distribution $\Delta$.   Then every other leaf
touched by this geodesic is also crossed
perpendicularly \cite{B-R4, Rampazzo2, Re1, Re2}.}
while the dynamical equations still
contain a term which is quadratic in $\dot u$.
This is the case of the Roller Racer,
our first example in Section~\ref{examples}.

Given a holonomic system  which is fit for jumps, 
Theorem \ref{fitprop} below provides a sufficient condition  in order that the  system  
remains fit for jumps also
 after the addition of  non-holonomic constraints.   This happens in the case of the rolling  ball,
our second example
in Section \ref{examples}. }
\end{remark}

\begin{thm}\label{fitprop} Assume that the external force $F$ does not depend on 
velocity\footnote{We can relax this hypothesis by allowing $F$ to be affine in the velocity.}. 
Suppose that the following conditions are satisfied:
 \begin{itemize}
  \item[(i)] The orthogonal distribution
$\Delta^\perp$ has zero curvature, i.e.: without  the nonholonomic constraint $\dot q\in \Gamma$
the system would be fit for jumps.
 \item[(ii)] There exists a $\Delta$-adapted system of coordinates $(q)= (q_\sharp, q_\flat)$,
 with $q_\sharp\doteq (q^1,\dots,q^N)$,  $q_\flat\doteq (q^{N+1},\dots,q^{N+M})$ , $\Delta=span\{\frac{\partial}{\partial q^1},\dots,\frac{\partial}{\partial q^N} \}$,
  such that the Hamiltonian $H$ is independent of the coordinates $q_\flat$.
     \item[(iii)] There
exists a basis $\{{\bf V}_1,\dots,{\bf V}_{N+M}\}$ for $T\Q$ such that, in this basis,
the projection ${\mathcal P}_I^*$ is represented by a
constant (i.e., independent of $\q$) matrix $P_I^*$.
\end{itemize}
Then,  even after the addition of
nonholonomic constraints $\dot q\in\Gamma$, 
the system remains fit for jumps .\end{thm}

{\it Proof.} We need to show that the time derivative
$\dot p_I$ is an affine function of $\dot u$.  Recalling (\ref{pIdot}) we have
\bel{sum}
\dot p_I~=~\left({\partial P_I^*\over\partial q}
\cdot \dot q\right) (p)
+ P_I^*(\dot p)  \,.
\eeq
By the assumption (iii), ${\partial P_I^*\over\partial q}=0$. Therefore
$$\dot p_I~=~P_I^*(\dot p),$$
 where $p(\cdot)$ solves the second equation in (\ref{Hamno}), namely
$$
\dot p ~= ~-\frac{\partial H}{\partial q} + F +R.
$$
By the hypotheses on $F$ and since $R_I=0$ it is then sufficient to show that 
 $\ds-\frac{\partial H}{\partial q} $ can be expressed as a  function of $(q,p_I, \dot u)$ 
 which is affine in $\dot u$.
 The assumption  (ii) yields
  $$
 -\frac{\partial H}{\partial q}~ =~ \left(-\frac{\partial H}{\partial q_\sharp}, -\frac{\partial H}{\partial q_\flat}\right)~ =~  \left(-\frac{\partial H}{\partial q_\sharp}, 0\right).
 $$
 On the other hand, by considering the cotangent natural basis $(dq_\sharp, dq_\flat)$ and the corresponding adjoint coordinates $(p_\sharp, p_\flat)$, we can express $\ds -\frac{\partial H}{\partial q_\sharp}$ as a function of $(q, p_\sharp,\dot q_\flat)$,
 $$
 -\frac{\partial H}{\partial q_\sharp} ~=~ -\frac{\partial H}{\partial q_\sharp}\Big(q,\,p_\sharp,
 \, p_\flat(q, p_\sharp,\dot q_\flat)\Big),
 $$
 where  $p_\flat  = p_\flat(q, p_\sharp,\dot q_\flat)$ is the function obtained by partially inverting the linear isomorphism $\ds \dot q 
 = \frac{\partial H}{\partial p}$. By the assumption (i), 
 the function  $\ds -\frac{\partial H}{\partial q_\sharp} $  (considered as a function of the 
 variables  $(q,p_\sharp, \dot q_\flat)$) is in fact affine w.r.t.~$\dot q_\flat$ (see e.g. \cite{B1} or \cite{Rampazzo2}). This is equivalent to saying that it is affine w.r.t.~$\dot u$, 
 because the map  $\dot u\mapsto \dot q_\flat(\dot u)$ is a linear 
 isomorphism\footnote{For a suitable choice of local  coordinates on $\U$,  this isomorphism coincides with the identity.}.   
 Let $\pi_\sharp$ denote the projection $p\mapsto \pi_\sharp(p) = (p_\sharp,0)$. By (\ref{c23}) we have $$(p_\sharp,0) ~=~ \pi_\sharp(p) ~=~ \pi_\sharp (p_I +p_{III})~ =~ \pi_\sharp (p_I + k(\dot u)).$$  
 Therefore
$-\frac{\partial H}{\partial q_\sharp} $  is  a function of $(q,p_I, \dot u)$ which is affine in  $\dot u$.
\endproof

\vskip0.51truecm
\section{A geometric interpretation of the quadratic term}
\setcounter{equation}{0}

We begin this section by observing that
it is not restrictive to assume that the manifold $\Q$
is an open subset of an Euclidean space.

Indeed, by the isometric embedding theorem
for Riemann manifolds
\cite{Cartan, Nash}, there exists a positive integer $M'$ and an embedding
$\iota: \Q\mapsto \R^{N+M+M'}$ (with $\R^{N+M+M'}$ endowed with the usual Euclidean metric) which preserves the Riemann metric on $\Q$  .
Let $\{ \bfe_1,\ldots,\bfe_{N+M+M'}\}$ be the standard orthonormal
basis of $\R^{N+M+M'}$.  Fix any point $\bar\q\in\Q$.
By the implicit function theorem, we can assume that
in a neighborhood of $\iota(\bar \q)$
the image of $\Q$ can be written as
$$\Big\{ (x_1,\ldots, x_{N+M+M'})\,;~~~x_i =
f_i(x_1,\ldots, x_{N+M})\,,~~~N+M+1\leq i \leq
N+M+M'\Big\}.$$
For $\bfx=(x_1,\ldots, x_{M+N+M'})$, we denote by $\bfq(\bfx)$
the unique point $\bfq$
in a neighborhood of
$\bar\bfq\in\Q$
such that the first $N+M$ coordinates of $\iota(\bfq)$ coincide with
 $(x_1,\ldots, x_{N+M})$.

 Next, consider the enlarged control manifold $\U'\doteq \U\times \R^{M'}$.
Together with the submersion $\pi:\Q\mapsto\U$,
consider  the submersion $\pi':\R^{N+M+M'}\mapsto \U'$
by setting
\bel{pi'd}\pi'(x_1,\ldots, x_{N+M+M'}) ~\doteq~\Big(\bfu,\,
 (u_{M+1}, \ldots, u_{M+M'})\Big)~\in~
\U\times \R^{M'},\eeq
where
$$\bfu = \pi(\bfq(\bfx)),\qquad\qquad
u_i~=~x_i - f_i(x_1,\ldots, x_{N+M})\quad\hbox{for}\quad N+M+1\leq i \leq
N+M+M'.$$
Calling $D\iota$ the differential of the embedding $\iota$ at the point $\bfq(\bfx)$,
we use  $\Delta'_\bfx \doteq D\iota\circ \Delta_{\bfq(\bfx)}$
to denote the holonomic distribution on $\R^{N+M+M'}$
generated by the submersion (\ref{pi'd}).
Given the non-holonomic distribution $\Gamma$ on $\Q$,
the corresponding non-holonomic distribution $\Gamma'$
 on $\R^{N+M+M'}$
is defined as
\bel{nonh'}\Gamma'_\bfx ~\doteq~\hbox{span}\Big\{ D\iota\circ \Gamma_{\bfq(\bfx)}\,,
~\bfe_{N+M+1},~\ldots~,\bfe_{N+M+M'}\Big\}.\eeq
Observe that, if (\ref{transv}) holds,
then we also have
$$\Delta'_\bfx +\Gamma'_\bfx~=~\R^{N+M+M'}.$$
Trajectories of the non-holonomic system with active constraints defined at
(\ref{nhol})-(\ref{submersion}) can now be recovered as trajectories
of the non-holonomic system with active constraints (\ref{pi'd})-(\ref{nonh'})
corresponding to control functions
$$t~\mapsto ~\Big(\bfu(t),\, u_{M+1}(t),\ldots,
u_{M+M'}(t)\Big)~\in ~ \U\times \R^{M'}$$
with $u_{M+1}(t)=\cdots=
u_{M+M'}(t) = 0$ for every time $t$.
\v
Based on the previous remarks, without loss of generality
we now analyze the geometric
meaning of the term $\Psi[\dot u, \dot u]$ in (\ref{centrifugal}),
assuming that $\Q=\R^{N+M}$ is a finite dimensional
vector space with Euclidean metric.
In this case the coefficients of the metric are constant:
$\bfg_{rs}\equiv g^{rs} \equiv  \delta_{rs}$. Moreover, with canonical
identifications we have
\bel{sidd}
p ~=~\dot q\,,\qquad\quad k(\dot u) ~=~h(\dot u)\,,
\qquad\quad
{\partial g^{-1}\over\partial q}~\equiv ~0\,\qquad\quad P_I=P^*_I\,.\eeq
Therefore, (\ref{thIdef})-(\ref{centrifugal}) yield the simpler
formula
\bel{Psid}
\Psi[\dot u, \dot u]~=~\left({\partial P_I\over\partial q}\cdot
h(\dot u)\right) (h(\dot u))\,.\eeq

In this setting, consider a system of local coordinates $(u_1,\ldots, u_M)$
on $\U$. To fix the ideas, let ${v} = (1,0\dots,0)\in\R^M$.
We seek a geometric characterization of the term
\bel{Psid2}
\Psi[{v}, {v}]~=~\left({\partial P_I\over\partial q}\cdot
h({v})\right) (h({v}))\,.\eeq
At each point $q\in\R^{M+N}$ in a neighborhood of $\bar q$,
choose an orthogonal basis of unit vectors $\V_1(q),\ldots, \V_{N+M}(q)$ as in
(\ref{familyVF}), but with $\V_{N+1}$ being the unit vector parallel to $h({v})$,
so that
$$h({v}, q) ~=~ \xi(q) \V_{N+1}(q)$$
for some $\xi(q)>0$.  Notice that we now write $h=h({v}, q)$, since we regard
$h$ as a vector field on $\Q=\R^{N+M}$.
In the following, we use $\langle\cdot,\cdot\rangle$ to denote the Euclidean inner product, and we write
$ (D\phi)\, \V_i$
for the directional
derivative of $\phi$ in the direction of $\V_i$.
Moreover, the map
$t\mapsto (\exp t \V_i)(\bar x)$ denotes the solution of the Cauchy problem
$$\dot x(t)~=~\V_i(x(t))\,,\qquad\qquad x(0)~= ~\bar x\,.$$
Computing (\ref{Psid2}) in terms of the basis $\{ \V_1,\ldots, \V_{N+M}\}$
we obtain
$$
P_I(w) ~=~\sum_{i=1}^{N-\nu} \la w\,,~\V_i\ra\, \V_i\qquad \forall w\in \R^{N+M}, $$
\bel{Psi3}\bega{rl}
\Psi[{v}, {v}]&\ds=~\sum_{i=1}^{N-\nu} \la \xi \V_{N+1}\,,(D\V_i)~\xi \V_{N+1}
\ra\, \V_i + \sum_{i=1}^{N-\nu} \la \xi \V_{N+1}\,,~\V_i\ra
\,(D\V_i)\xi \V_{N+1}
\cr
&\ds=~\sum_{i=1}^{N-\nu} \la \xi \V_{N+1}\,,~(D\V_i)\xi \V_{N+1}
\ra\, \V_i \,.\enda\eeq
We now claim that, for $1\leq i \leq N-\nu$, the Lie bracket
$$[\xi\V_{N+1}\,,~\V_i] ~=~(D\V_i)\xi \V_{N+1} -(D(\xi \V_{N+1}))\V_1$$
lies in the subspace $\Delta_\bfq$ and hence is orthogonal to $\V_{N+1}$.
Indeed, for every $\bar q$ one has
$$[\xi\V_{N+1}\,,~\V_i](\bar q)~=~\lim_{\ve\to 0}~{q(\ve)-\bar q\over\ve^2}\,,$$
where
$$q(\ve)~\doteq~
\Big(\exp (-\ve\V_i)\Big)\circ\Big(\exp (-\ve\xi\V_{N+1})\Big)\circ
\Big(\exp (\ve\V_i)\Big)\circ\Big(\exp (\ve\xi\V_{N+1})\Big)(\bar q)\, .$$
However, for any $q$ we have
$${d\over dt} \pi \Big(\exp (t\V_i)\Big)(q) ~=~0\,,
\qquad\qquad {d\over dt} \pi\Big(\exp (t\xi\V_{N+1})\Big) (q) ~=~{v}\,.$$
Hence, for every $\ve$ one has
$$\pi(q(\ve))~=~\Big(\exp (-\ve {v})\Big)\circ
\Big(\exp (\ve {v})\Big)(\pi(\bar q)) ~=~\pi(\bar q)\,.$$
This implies
$$[\xi\V_{N+1}\,,~\V_i](\bar q) ~\in~ ker\,D\pi (\bar q)~=~\Delta_{\bar q}\,,$$
as claimed. Therefore
\bel{Lieb}\la \V_{N+1}\,,(D\V_i)(\xi \V_{N+1})
\ra~=~ \la \V_{N+1}\,,~(D(\xi \V_{N+1}))\V_i\ra\,.\eeq
Next, we observe that
$$\la \V_{N+1}\,,~(D\V_{N+1})\V_i\ra ~=~{1\over 2} (D
|\V_{N+1}|^2)\V_i ~=~0\,,$$ because all vectors $\V_{N+1}(q)$ have unit length.
{}From  (\ref{Psi3}) and (\ref{Lieb}) we obtain
\bel{Psi4}\bega{rl}
\Psi[{v}, {v}]
&\ds=~\sum_{i=1}^{N-\nu} \la \xi \V_{N+1}\,,~(D (\xi \V_{N+1}))\V_i
\ra\, \V_i \cr
&\cr &=\ds ~\xi\, \sum_{i=1}^{N-\nu} ((D\xi)\V_i) \, \V_i
~=~\xi \, P_I^*(\nabla \xi)\,.\enda\eeq

Recalling (\ref{argmin}),
we can think of $\xi(q)= |h({v},q)|$ as the infinitesimal
distance between  leaves of the foliation generated by the submersion
$\pi$ at the point $q$, computed along paths that respect the
non-holonomic distribution.
According to (\ref{Psi4}), if the directional derivative
of  this distance is zero  in all directions of the subspace
$(T_q\Q)_I=\Delta_q\cap\Gamma_q$, then
the term $\Psi[\bfv,\bfv ]$ vanishes for all $\bfv\in T_{\pi(q)}\U$.
\v
Going back to the general case of a Riemann manifold $\Q$ and a submersion
$\pi:\Q\mapsto\U$, we thus obtain
the following geometric characterization.
For every $\bfu\in \U$, consider the leaf
\bel{leaf1}
\Lambda_\bfu~\doteq~\{ \q\in\Q\,;~~\pi(\bfq) = \bfu\}\,.\eeq
Then for every $\bfq\in\Lambda_\bfu$ and every
tangent vector $\bfv\in T_\bfu\U$ there exists a unique
vector $\bfh(\bfv)\in (T_\bfq\Q)_{III}$ such that $D\pi\cdot \bfh(\bfv)=\bfv$.
Recalling (\ref{argmin}),
one could define an inner product on $\langle\cdot,\cdot\rangle_\bfu$
on $T_\bfu\U$ by setting
\bel{ip6}
\langle \bfv,\bfv\rangle_\bfu~\doteq~\bfg_\bfq[ \bfh(\bfv),\, \bfh(\bfv)]\,.
\eeq
Of course, in general this inner product depends on the choice of
the point $\bfq$ along the leaf $\Lambda_\bfu$.
According to (\ref{Psi4}), the condition $\Psi\equiv 0$
means that the directional derivative of the right hand side of
(\ref{ip6}) vanishes along $\Delta_\bfq\cap\Gamma_\bfq$.  More precisely:
\v
\begin{prop}
The following conditions are equivalent.
\begi
\item[(i)]
$\Psi\equiv 0$.

\item[(ii)] For every $\bfu\in \U$, $\bfv\in T_\bfu\U$, and
every smooth path  $s\mapsto \bfq(s)\in \Lambda_\bfu$
such that $\dot \bfq(s)\in\Delta_{\bfq(s)}\cap \Gamma_{\bfq(s)}$, the inner
product  $\bfg_{\bfq(s)}[ \bfh(\bfv),\, \bfh(\bfv)]$ is constant.
\endi
\end{prop}
\v

\begin{remark} {\rm In general, the condition  $\Psi\equiv 0$ does not guarantee
the existence of a Riemannian metric on $\U$ defined by (\ref{ip6}).
In order that this metric be well defined, i.e.~independent of the
point $\bfq\in \Lambda_\bfu$, one also needs to assume that
the distribution
$\Delta_\bfq\cap\Gamma_\bfq$ is completely non-integrable restricted to
each leaf $\Lambda_\bfu$.  More precisely, every two points $\bfq_1,
\bfq_2\in \Lambda_\bfu$
should be connected by a path $t\mapsto \bfq(t)\in \Lambda_\bfu$ with $\dot \bfq(t)\in
\Delta_{\bfq(t)}\cap\Gamma_{\bfq(t)}$ for all $t$.
Notice that this condition is clearly satisfied if
the non-holonomic constraint is not present
(i.e., ~$\Gamma_\bfq = T_\bfq\Q$) and the leaves $\Lambda_\bfu$ are
connected.  From (\ref{Psi4}) we thus recover some earlier characterizations
of the property ``fit for jumps" given in  \cite{AB1,  Rampazzo2}.}
\end{remark}

\section{Two examples}\label{examples}
\setcounter{equation}{0}

\subsection{The Roller Racer}

As a  toy model, consider the {\it Roller Racer},
which will serve as a simple illustration of the  general theory.
This  is a classical example of a non-holonomic system,
widely  investigated within the  theory of the
{\it momentum map} \cite{Bloch, BL}.
It consists of two rigid planar bodies,
connected at a point by a rotating joint,
as shown in Figure 2. One of the two bodies has a
much larger mass than the other.
Let $\rho$ be the distance between
the joint and the center of mass of the heavier body.
To simplify computations we also assume that the center
of mass of the lighter body coincides with the joint.

The coordinates $(q^1,q^2,q^3,u)$ are as follows .
We let $(q^3,q^1)=(x,y)$
be the Euclidean coordinates of the center of mass of the large body.
Moreover, $q^2=\theta$ is the counter-clockwise
angle between the horizontal
axis and the major axis of the large body, and
$u=\phi$ is the counter-clockwise angle between the
 major axes of the two  bodies.

 The non-holonomic constraint consists in assuming that each
 pair of wheels roll without slipping,
 parallel to the corresponding body.
This corresponds to the condition
 $$(\dot q, \dot u)~\in~ \Gamma ~\doteq~ ker(\om^1)\cap ker (\om^2),
 $$
 where
 \bel{om12}\left\{\begin{array}{l}
 \om^1\doteq \cos q^2\,dq_1-\sin q^2 \,dq^3\,,
 \cr
 \om^2\doteq \cos(q^2+u)\,dq_1+ \rho\,\cos u\, dq^2
 -\sin(q^2+u)\, dq^3
\,.\end{array}\right.
\eeq
In this case the (non integrable) distribution is 2-dimensional.
The admissible motions $t\mapsto(q,u)(t)$ are thus subject to
\bel{RRc}\left\{
 \begin{array}{l}
  \cos q^2(t) \,\dot q^1(t)  -\sin q^2(t)\, \dot q^3(t)~= ~0\,,
 \cr\,\cr
\cos(q^2(t)+u(t))\,\dot q^1(t)+ \rho\,\cos u(t)\,
\dot q^2(t) -\sin(q^2(t)+u(t))\, \dot q^3 (t)~= ~0\,.
\end{array}\right.
\eeq
\begin{figure}
\centering
\includegraphics[scale=0.50]{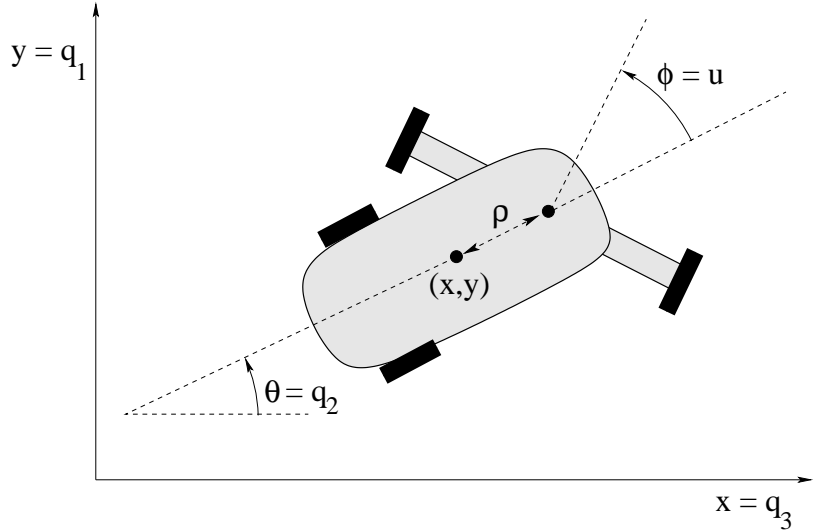}
\label{f:a91}
\caption{The roller racer.}
\end{figure}

The coordinate $\phi= u$ is regarded as a control. Notice that
the transversality condition (\ref{transv}) is
trivially  satisfied, because
$$
\Delta^{ker}~=~span\{du\}, \qquad
\Gamma^{ker}~=~span\{\om^1,\om^2\},
$$
and $span \{\om^1,\om^2\}\cap span \{du\}= \{0\}$.

For simplicity, we assume that
the large body has unit mass.
The  moment of inertia of the large body w.r.t.~the
vertical axis passing through
 its center of mass
is denoted by $I$.
The mass of the small body is regarded as negligible, but its
moment of inertia $J$ w.r.t.~the vertical
axis passing through the center of mass (of the small body)
is assumed to be different from zero.\footnote{This is a standard
approximation adopted in the existing literature.}
The  kinetic matrix ${\bf G}=(g_{i,j})$ and its inverse ${\bf G^{-1}}=(g^{i,j})$ are computed as
$$
(g_{i,j})~ =~\pmatrix{ 1& & 0&&0 &&0\cr
\cr
0 &&I+J &&0&&J\cr
\cr
0 &&0 &&1&&0
\cr\cr 0
&&J &&0&&J}  \qquad  (g^{i,j}) ~=~\pmatrix{1&&0&&0&&0\cr\cr
0&&{1\over I}&&0&&-\frac{1}{I}
\cr\cr 0 &&0& &1&&0
\cr\cr 0 &&-\frac{1}{I} &&0&&\frac{I+J}{IJ}}.
$$

Let us start by finding a basis for the decomposition $T\Q=(T\Q)_{I}\oplus (T\Q)_{II}\oplus (T\Q)_{III}
$. Notice that the vectors
 $$\bfw_1=  2\rho\cos{u} \sin{q^2}~\frac{\partial}{\partial q^1} +~ 2 \sin{u}\frac{\partial}{\partial q^2} +~ 2\rho \cos{u} \cos{q^2} \frac{\partial}{\partial q^3}\,, ~\qquad~~
\bfw_2= \frac{\partial}{\partial u}\,,
$$
form  a basis for $\Gamma$. Since $\Delta
 = span\{\frac{\partial}{\partial q^1}, \frac{\partial}{\partial q^2}, \frac{\partial}{\partial q^3}\}$,
$\{\bfw
_1\}$ is a basis for $(T_\Q)_I =\Gamma\cap\Delta$.

It is also straightforward to verify that
\bel{gammaperp}
(T\Q)_{II}~ = ~\Gamma^{\perp} ~=~ \ker\left\{\g(\bfw_1),\g(\bfw_2)\right\} ~=~ span\left\{ \bfv_2, \bf
v_3\right\}\,,
\eeq
where
$$ \bfv_2~=~{I \csc{q^2} \tan{u}\over \rho}~ \frac{\partial}{\partial q^1} - \frac{\partial}{\partial q^2} + \frac{\partial}{\partial u},\qquad \bfv_3~=~-\cot{q^2}~ \frac{\partial}{\partial q^1}+ \frac{\partial}{\partial q^3}\,.$$

Since $(T\Q)_{III}$ is orthogonal to $(T\Q)_{I}\oplus (T\Q)_{II}$,
\bel{III}
(T\Q)_{III}~=~\ker\Big\{\g({\bfw_1}), \g(\bfv_2),
\g(\bfv_3 )\Big\}=~span\{\bfv_4\}
\eeq
with
$$\bv_4~\doteq~{1\over \Delta_0}\left( -{1\over 2}J\rho\sin{q^2}\sin{2u}~\frac{\partial}{\partial q^1}~-J\sin^2{u}~\frac{\partial}{\partial q^2}~-{1\over 2}J\rho\cos{q^2}\sin{2u}~\frac{\partial}{\partial q^3}\right)+\frac{\partial}{\partial u}\,,
$$
we have
 $\Delta_0 \doteq \rho^2\cos^2{u}+(I+J)\sin^2{u}$.

Defining $\V_1\doteq \bfw_1, \V_2\doteq\bfv_2, \V_3\doteq\bfv_3,
\V_4\doteq \bfv_4$ we obtain a family of vector fields (on an open subset of the configuration manifold) as in (\ref{familyVF}).

To obtain the equations of motion
 we need to compute the coefficients (\ref{thetatilde}).
 For this purpose, let us begin with observing that the form ${\bf \Omega}_1$
 generating $(T^*\Q)_I$ is given by
$$\begin{array}{l}{\bf \Omega}_1~ = ~\Omega_{1,1}dq^1 +
\Omega_{2,1}dq^2+
\Omega_{3,1}dq^3+
\Omega_{4,1}dq^4 ~ =~ \g(\V_1)
\\\,\\
\qquad =~(2\rho\cos{u} \sin{q^2})dq^1 +
(2(I+J) \sin{u})dq^2+
(2\rho \cos{u} \cos{q^2})dq^3+
(2J \sin{u})dq^4.
\end{array}
$$
Let us recall that the projection matrix $P^*_I$ is computed as (see (\ref{copromat}))
$$
{(P_{I}^*)^r}_s = \sum_{k=1}^{4}\frac{g^{r,k}\Omega_{k,1}
 \Omega_{s,1}}{\sum_{a,b=1}^{4}g_{a,b}{V^a}_1
 {V^b}_1}\,.
 $$
Let us set $\Delta_1= I+J+\rho^2+(-I-J+\rho^2)\cos{2u}$. Then an elementary
computation of (\ref{Hamclosedx1}) yields  the following system of
four differential equations:
\bel{RollerRacer}
\left\{\begin{array}{l}
\dot q^1~=~2\rho\cos{u}\sin{q^2}\cdot \xi -{J\rho\sin{q^2}\sin{2u}\over 2\Delta_0}\cdot \dot u\,,\\\,
\\
\dot q^2~=~2\sin{u}\cdot \xi -{J\sin^2{u}\over\Delta_0} \cdot \dot u\,,\\\,
\\
\dot q^3~=~2\rho\cos{q^2}\cos{u}\cdot \xi -{J\rho\cos{q^2}\sin{2u}\over 2\Delta_0} \cdot\dot u\,,\\\,
\\
\dot \xi~=~-{2(I+J-\rho^2)\sin{2u}\over \Delta_1}\cdot \xi\,\dot u+{8J\rho^2\cos{u}\over\Delta_1^2}\cdot \dot u^2.
\end{array}\right.
\eeq

\begin{remark}{\rm
Assume that we implement a vibrating control, of the form $u(t) = \bar u+\ve K\sin(t/\ve)$
for $K\in\R$ some $\ve>0$ small.  By the results from \cite{B-R3, B-R4},
each solution of the system
\beq \left\{\begin{array}{l}\ds
\dot q^1~=2\rho \cos{\bar u} \sin q^2\cdot\xi
\\\,\\
\dot q^2~=2 \sin\bar u\cdot \xi
\\\,\\
\dot q^3~=2\rho \cos{\bar u} \cos q^2\cdot\xi
\\\,\\
\ds\dot \xi~=~ {4J\rho^2\cos{\bar u} \cdot K^2\over\Delta_1^2}
\end{array}\right.
\eeq
(formally obtained by neglecting higher order terms and then averaging),
can be uniformly approximated by the solutions of the original system
\ref{RollerRacer}
\footnote{More generally, one could approximate any solution of a (possibly) impulsive system obtained  \ref{RollerRacer} by replacing $(\dot u,\dot u^2)$  with a pair $(a,\mu)$, where $a$ is a $L^1$ control and $ \mu$ is a positive measure such that $a^2\leq\mu$ (in the sense of measures)}.
It is then easy to verify a kinematical well-known behavior of the Roller Racer: fast oscillations of the handlebar around a angle $\bar u$ produce  a motion which, asymptotically, is a superposition of forward motion  and a rotation of the angle $q^2$ (which is the angle between the principal axis of the large body and the $x$-axis, namely the $q^3$ axis).   In particular if the oscillations are around the central position $u=0$, then $q^2$ is constant and the motion approaches asymptotically a forward motion.
Of course this agrees with results obtained with methods based on momentum maps, where, in particular,   controls are forces or torques.  See
\cite{Bloch} or the extensive analysis of the Roller Racer  in \cite{KriTsa}.}
\end{remark}

\subsection{A ball rolling on rotating disc}

The next example illustrates an application of Theorem \ref{fitprop}.
It represents a modification of the classical case (see~\cite{Bloch}) of a
mechanical system consisting of a ball that rolls without sliding on a disc which rotates with constant speed (Figure 3). The somehow surprising fact is that the
ball does not move away toward infinity, regardless of the angular speed
of the disc. In fact, the barycenter of the ball remains confined within a bounded region, depending on the initial conditions.
The ``surprise" actually comes from the fact that intuition suggests a sort of centrifugal effect, pushing the  ball outward in the radial
direction.

We slightly modify this example by allowing  the rotational speed of the disc
to be {\it not constant}. More precisely,  we consider
the rotational angle $u$ of the disc as a (state-space) parameter
$u$, regarded as a control. We will show that, as in the case of constant rotational speed, a small rapid oscillation of the disc around the center does not
instantly push the ball toward infinity.
As a consequence, the system is ``fit for jumps".
 To recast this system in our general framework,
consider a material disc $D$, centered at a point $O$  with inertial moment $I>0$.
If $u(t)$ designates the counter-clockwise
angle of the disc with respect to inertial coordinate axes,
its kinetic energy is given by
$\frac{1}{2} I\dot u^2$.  Actually, in our setting the value of $I$ is
irrelevant for the motion of the system, because the angle $u$
is taken as a control variable.
The configuration manifold  of the rotating disc together with the rolling ball
can be identified with $\Q\doteq \R^2\times SO(3)\times S^1$.
The first two coordinates $(x,y)\in\R^2$
describe the position of the contact point
w.r.t~the inertial frame.  A unitary matrix $R\in SO(3)$ represents the rotation
of the  ball, while the angle $u\in S^1$ describes the rotation of the disc.
 Here the submersion describing the control constraint is simply the projection
 $$
 \pi: \Q~\to ~ S^1,\qquad\qquad \pi(x,y,R,u)~=~ u\,.$$
In addition, the non-holonomic constraints are represented by the linear equations
 \bel{nonhol-sph}\left\{\begin{array}{l}
\dot x + r\omega_2 + \dot u y~=~ 0\,,  \\\,\\
\dot y -r\omega_1 - \dot u x ~=~ 0\,,
\end{array}\right.
\eeq
where $r$ is the radius of the  ball and $\omega=
(\omega_1,\omega_2,\omega_3)$ is the angular velocity vector of the
ball  w.r.t.~the inertial frame. For simplicity, we assume that the
inertial moment of the disc w.r.t.~its center is $I=1$ and
 that the ball has homogeneous density and total mass equal to $1$.
The kinetic energy ${\it T}$ of the  system is then given by
\bel{ke4}
{\it T} ~= ~\frac{1}{2} \left(\dot x^2+\dot y^2+\kappa^2
\sum_{i=1}^{3}\omega_i^2 + \dot u^2\right),
\eeq
where $\kappa^2$ is the moment of inertia of the ball w.r.t.~any of its axes.

\begin{figure}
\centering
\includegraphics[scale=0.50]{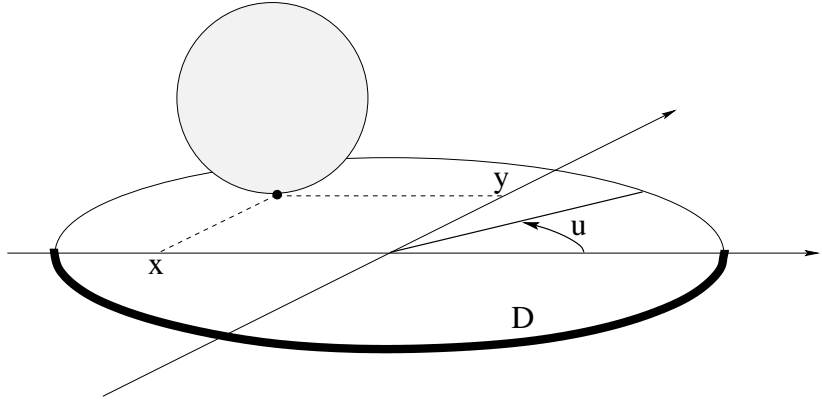}
\label{f:a92}
\caption{A ball rolling without sliding on a rotating disc.}
\end{figure}

Let $(q^1,\dots,q^6)$ be local coordinates for $\Q$, where $q^4\doteq x$,
$q^5\doteq y$, $q^6\doteq u$, while $(q^1,q^2,q^3)$ are local coordinates
for $SO(3)$. Recalling the meaning of angular velocity $\omega$
and the form of the kinetic energy $\it T$,  we can find
an orthonormal basis $\{\V_1,\ldots, \V_6\}$
of the tangent space $T\Q$ such that the following holds:
$$\V_i~\doteq~\frac{\partial}{\partial q^i}\qquad\qquad \hbox{for}~~~i\in \{4,5,6\}\,,$$
\bel{gliA}
\V_j~=~\sum_{h=1}^3 A^h_j(q^1,q^2,q^3)\frac{\partial}{\partial q^h}\qquad
\qquad \hbox{for}~~~j\in\{1,2,3\}.
\eeq
Moreover, if $\omega = (\omega_1,\omega_2,\omega_3)$
describes the angular velocity of the ball, then
$$
(\dot q_1, \dot q_2,\dot q_3)  ~=~ \sum_{j=1}^3\kappa\,\omega_j \V_j\,.
$$
In terms of these coordinates, if $\dot q = \sum_j c_j \V_j$, the non-holonomic constraints (\ref{nonhol-sph})
take the form
$$
c_4 + {r\over \kappa} c_2 + q^5 c_6~=~0\,,\qquad\qquad
c_5 - {r\over\kappa} c_1 - q^4 c_6~=~0\,.$$
We thus obtain
 \bel{dec8}
  \Delta~=~span\{\V_1,\V_2,\V_3,\V_4,\V_5\},\qquad\qquad \Delta^\perp~ =~span\{\V_6\}, \eeq
\bel{dec9}  \Delta\cap \Gamma ~= ~span\left\{\V_1 +{r\over \kappa} \V_5,
   \, \V_2- {r\over\kappa} \V_4,\, \V_3\right\},\eeq
  \bel{10}  (\Delta\cap \Gamma)^\perp ~=~ span\left\{\V_1 - {\kappa\over r} \V_5,
    \, \V_2+ {\kappa\over r} \V_4,\, \V_6\right\}.
    \eeq
 Using the basis$\{\V_1,\dots,\V_6\}$, all projections
    $P_J, P^*_J$, $J\in \{I,II,III\}$ are given by constant matrices, so hypothesis (iii) in Theorem \ref{fitprop} is verified.
By (\ref{ke4})-(\ref{dec8}) it is clear that the matrix $g$ of the kinetic energy with respect to the basis $(\frac{\partial}{\partial q^1},\dots,\frac{\partial}{\partial q^6})$ is independent of $q^6$ and has all zeros in the sixth row and the six column, except for the corner entry which is equal to $1$. Therefore the same is valid for the inverse matrix $g^{-1}$. 
Denoting by $p^T$  the transpose of the row vector $p$, this implies  
\begi
\item  $H \doteq \frac{1}{2} pg^{-1}p^T$ is independent of $q^6$.

\item If the
non-holonomic constraint  (\ref{nonhol-sph}) is not present, then the  system is fit for jumps. 
\endi

By Theorem \ref{fitprop} we thus conclude that this system is fit for jumps. 
In particular, a rapid small oscillation of the disc
around its center will not cause the ball to fly away to infinity.

\vs

{\bf Acknowledgements.}
The research of the first author was
supported by NSF   through grant DMS-1108702, ``Problems of Nonlinear Control".  The third author was
 partially supported by the Italian Ministry of  University and Research, PRIN:  ``Metodi di viscosit\`a, di teoria del
controllo, e di analisi nonsmooth per PDE nonlineari ellittico-paraboliche e di tipo Hamilton-Jacobi".

\end{document}